\documentclass[12pt]{amsart}
 \usepackage{amssymb,amsfonts,amsmath,amsopn,amstext,amscd,latexsym, amsthm, enumerate,mathrsfs}
 \usepackage{float}
 \usepackage{color}
 \usepackage{xcolor}
 \usepackage{hyperref}
 \hypersetup{
 	colorlinks   = true,
 	citecolor    = blue
 }
 
 \usepackage[margin=30mm]{geometry}
 \usepackage{bbm}
 \usepackage[all,cmtip]{xy}
 \usepackage{graphicx}
 \usepackage{enumerate}
 \usepackage[shortlabels]{enumitem}
 
 \newtheorem{theorem}{Theorem}[section]
 
 \newtheorem{lemma}[theorem]{Lemma}

 \newtheorem{lem}[theorem]{Lemma}
 
 \newtheorem{proposition}[theorem]{Proposition}

 \newtheorem{cor}[theorem]{Corollary}
 \newtheorem{definition}[theorem]{Definition}
 
 \theoremstyle{remark}
 
 \newtheorem{remark}[theorem]{Remark}

 \numberwithin{equation}{section}

\begin{document}
 	
 	\title{Characterization of $PSL(2,11)$ by the Set $U(G)$}

 	\author[M.~Hemmati]{Mina Hemmati Tirabadi}
 	\address{Department of Pure Mathematics, Faculty of Mathematical Sciences, Tarbiat Modares University, 14115-137, Tehran, Iran}
 	\email{mina$\_$hemmati@modares.ac.ir (mht.199372@gmail.com)}
 	
 	\author[A.~Iranmanesh]{Ali Iranmanesh}
 	\address{Department of Pure Mathematics, Faculty of Mathematical Sciences, Tarbiat Modares University, 14115-137, Tehran, Iran}
 	\email{iranmanesh@modares.ac.ir}

 	\subjclass[2010]{ 20E45, 20E32}
 	
 	\date{\today}

 	\keywords{Conjugacy class; Simple Group; Conjugate type vector; Same-size \linebreak conjugacy class}

 	\begin{abstract}
In this paper, we prove that if $G$ is a finite simple group with the \linebreak same-size conjugacy class set $U(G) = U(PSL(2,11))$, then $G$ is isomorphic to $PSL(2,11)$.
	\end{abstract}
	\maketitle
	\section{Introduction}
	Let $G$ be a finite group and $n_i$ denote the distinct sizes of conjugacy classes of $G$, where $i \in \mathbb{N}$. In 2019, M.Zarrin and N.Ahmadkhah \cite{Zarrin} introduced a new equivalence relation on $G$ based on the sizes of its conjugacy classes. They focused on the set of sizes of the equivalence classes with respect to
	this relation, defining it as follows:
	For all $g,h \in G $, we write $ g \sim h $ if and only if $\mid g^G\mid = \mid h^G \mid$, where $\mid g^G \mid$ denotes the size of the conjugacy class of $g$.
The equivalence classes under this relation are $\{x \in G \quad ; \quad \mid x^G \mid = n_i\}$, and the sizes of these equivalence classes are denoted by $u_G(n_i)$. The set of all distinct values of $u_G(n_i)$ is denoted by $U(G)$, and is called the same-size conjugate class set of $G$. (In this paper, we refer to this set as the same-size conjugacy class set.)

M.Zarrin and N.Ahmadkhah used the same-size conjugacy class set to characterize certain simple linear groups, including $PSL(3,3)$ and $PSL(2,q)$, where $q \in \{5,7,8,9,17\}$. They proved that if $G$ is a simple group with $U(G)=\{1,pq,4p,8q\}$, where $p$ and $q$ are prime numbers, then $p$ and $q$ are distinct odd primes and $G$ must be isomorphic to $A_5$.

In this paper, we extend their work by exploring the characterization of simple groups via their same-size conjugacy class set. We have proved, if $G$ is a finite simple group with $U(G) = U(PSL(2,11))$, then $G \cong PSL(2,11)$.
	\section{Preliminary}
	Throughout this paper, we assume $G$ to be a finite simple group. The set of prime divisors $\pi(G)$ is defined as the set of all prime divisors of the order of $G$. The groups, $A_n$, and $PSL(n,q)$ represent the Alternating group of order $n!/2$, the Projective \linebreak Special Linear group over a finite field of order $q$, respectively. Let $r$ denote the \linebreak number of distinct conjugacy classes of $G$, $x^G$ the conjugacy class of $x$, and $C_G(x)$ the centralizer of $x$. The set of all distinct  centralizers of $G$ is denoted by $Cent(G)$. The size of the conjugacy class $x^G$ is denoted by $n_i$, where $1 \leqslant i \leqslant r$. The set of \linebreak all conjugacy class sizes of $G$ denoted by $cs(G)$. If the conjugacy class sizes of $G$ are ordered in ascending order, the conjugate vector type and the conjugate type rank $G$, are defined as $V(G)=(n_1,...,n_r)$, where $n_i < n_j$ for $i<j$, and $r-1$, respectively. For \linebreak $u_G(n_i)=\mid \{ x \in G \quad ; \quad \mid x^G \mid = n_i\}\mid $, we define $U(G)=\{u_G(n_i) \quad ; \quad i \in \{1,...,r\}\}$. Clearly, the order of $G$ satisfies :
	\begin{align*}
		\mid G \mid = \sum_{i\in I} u_G(n_i).
	\end{align*}
	The following propositions and theorems, which are used throughout this paper, are fundamental to our discussion.
	\begin{proposition}\cite[Proposition 2.8]{Zarrin}\label{1.2}
		Let $G$ be a finite group and $n$ be a conjugacy class size of $G$. Then $n$ divides $u_G(n)$.\qquad $\Box$
	\end{proposition}
	Burnside’s Class Theorem is a central result relating conjugacy classes and group simplicity.
	\begin{theorem}\cite[Burnside's Class Theorem]{Burnside}\label{1.1}
If a finite group has a conjugacy class with prime power cardinality, the group is not simple.
    \end{theorem} 
    \begin{proposition}\cite[Lemma 3.7]{Zarrin}\label{1.4}
    	Let $Z(G)=1$. Then $\pi(G)=\bigcup_{i=1}^{r} \pi(n_i)$, where $V(G)=(n_1,n_2,...,n_r)$.
    \end{proposition}
    Some simple groups are characterized by their conjugate type rank. These \linebreak investigations are particularly useful in identifying specific group structures.
    \begin{theorem}\cite{Ito4}\label{1.3}
    	A simple group of conjugate type rank 4 is isomorphic to some $PSL(2,l)$, where $l$ is an odd prime power bigger than 5. The converse is also true.
    \end{theorem}
    Several authors have attempted to characterize certain simple groups, such as the Projective Special Linear groups $PSL(n,q)$, by their set of prime divisors and its size. The following results are particularly useful:   
   \begin{theorem}\cite[Theorem 4]{Herzog}\label{1.5}
   	Let $G\cong PSL(2,q)$, where $q = p^n >3$ and $p$ is a prime. Assume that $\mid G  \mid$ is divisible by three primes only. Then $G$ is isomorphic to one of the following groups:
   	\begin{enumerate}
   		\item [(i)]
   		$p = 2$ , $PSL(2,4)$ and $PSL(2,8)$;
   		\item[(ii)]
   		 $p > 2$ , $PSL(2,5)$ , $PSL(2,7)$ , $PSL(2,9)$ and $PSL(2, 17)$.
   	\end{enumerate}
     \end{theorem}
     \begin{definition}\cite{Bugeaud}\label{1.7}
     	The group $G$ is called $k_n$-group, if $G$ is simple group and $\mid \pi(G) \mid =n$, where $n \in \mathbb{N}$.
     \end{definition}
     \begin{theorem}\cite{Bugeaud}\label{1.6}
     	If $G$ is a simple $k_3$-group, then $G$ is isomorphic to one of the following groups:
     	
     		$A_5$, $A_6$, $PSL(2,7)$, $PSL(2,8)$, $PSL(2,17)$, $PSL(3,3)$,$U_3(3)$ and $U_4(2)$. 
     \end{theorem}
  \begin{theorem}\cite{Bugeaud}\label{1.10}
 	If $G$ is a simple $k_4$-group, then $G$ is isomorphic to one of the simple groups:
 	
 	$A_9$, $A_{10}$, $M_{11}$, $J_2$, $L_2(q)$ \quad where $q$ is a prime power \quad, $L_3(2^2)$\quad,$L_3(2^3)$\quad,$L_3(5)$\quad,$L_3(7)$\quad,$L_3(17)$\quad,$L_4(3)$\quad, $O_5(2^2)$, $O_5(3^2)$,$O_5(5)$,$O_5(7)$,$O_7(2)$,$O^{+}_8(2)$,$G_2(3)$,$U_3(2^2)$,$U_3(3^2)$,\linebreak $U_3(5)$,$U_3(7)$,$U_4(3)$,$U_5(2)$,$Sz(2^3)$,$Sz(2^5)$,$3D_4(2)$,$F_4(2)'$.
 	
 \end{theorem}
 \begin{theorem}\cite[Theorem 3.2]{Chen}\label{1.8}
 	Let $G$ be a group. $G \cong PSL(2,p)$ if and only if $\mid G \mid = \dfrac{p(p^{2}-1)}{(2,p-1)}$ and $G$ has a special conjugacy class size of $\dfrac{(p^{2}-1)}{(2,p-1)}$, where $p$ is a prime and not equal to seven.
 \end{theorem}
\section{main Theorem}
In this section, we examine the structure of a simple group $G$ where the set \linebreak $U(G) = \{1,rq,4rq,8pq,8pr\}$, with $p$, $q$, and $r$ representing prime numbers. Our \linebreak objective is to show that $p$, $q$, and $r$ must be distinct odd primes.
To achieve this, we will explore all possible relationships between the number $2$, $p$, $q$, and $r$, as outlined in the following lemmas. This analysis is crucial for understanding the connections under which the group $G$ can maintain its simplicity. 
\begin{remark}\label{2.1}
According to Theorem \ref{1.1} and Proposition \ref{1.2}, the following sets can not appear as the same-size conjugacy class sets for simple groups:
\begin{align*}
	 \{1 , r^2 , 4r^2 , 16r\}, \{1, p^2 , 4p^2 , 8p^2\}, and  \quad \{1 , r^2 , 4r^2 , 8pr\}.
\end{align*}
These sets do not satisfy the conditions required for conjugate sets in the context of simple groups, as established by the aforementioned results.
\end{remark}
\begin{lemma}\label{2.2}
	The following sets can not appear as the same-size conjugacy class sets for simple groups, where  $p$, $q$, and $r$, are distinct odd prime numbers: 
	\begin{align*}
	\{1,2p,8p,16p\}, \{1,2q,8pq,8q,16p\}, \{1,2p,8p,16p ,8p^2\}.
	\end{align*}
	These sets do not satisfy the necessary conditions for conjugate sets in the context of simple groups, as demonstrated by the constraints established in previous theorems.
	\end{lemma}
\begin{proof}
	Let set $U_1=\{1,2p,8p,16p\}$, $U_2= \{1,2q,8pq,8q,16p\}$ and $U_3=\{1,2p,8p,16p,8p^2\}$. By contradiction, let $G_i$ be a simple group corresponding to $U_i$, where $1 \leqslant i \leqslant 3$.\linebreak Suppose that:
	\begin{align*}
			u_{G_1}(n_2) = 2p, u_{G_2}(n_2) = 2q, \quad and \quad u_{G_3}(n_2) = 2p.
	\end{align*}
According to Theorem \ref{1.1} and Proposition \ref{1.2}, we have:
\begin{align*}
 n_2=2p\quad in \quad G_1, \quad n_2=2q \quad in \quad G_2 \quad and \quad n_2=2p \quad in \quad G_3.
\end{align*}
So $\mid G_1 \mid$, $\mid G_2 \mid$, and $\mid G_3 \mid$ are even numbers. On the other hand, the order of $G$, is given by $\mid G \mid = \sum_{i}(u_G(n_i))$, thus $\mid G_1 \mid$, $\mid G_2 \mid$, and $\mid G_3 \mid$ are odd numbers. This leads to a contradiction.
\end{proof}
Before we demonstrate that there is no finite simple group $G$ with \linebreak $U(G) = \{1,rq,16q,4rq,16r\}$, where $q$ and $r$ are distinct odd primes, we will first show that such this group must have 
$\mid cs(G) \mid = 5$.
\begin{lemma}\label{2.3}
	Let $G$ be a finite simple group with $U(G)=\{1,rq,16q,4qr,16r\}$, where $q$ and $r$ are distinct odd primes. Then, for any distinct indices $i$ and $j$, the size of the conjugacy classes of $G$ are different, i.e., $n_i \neq n_j$  for $i\neq j$, where that $n_i$ represents the size of the conjugacy class of $G$. 
\end{lemma}
\begin{proof}
Let $u_G(n_1)=1$, $u_G(n_2)=rq$, $u_G(n_3)=2q$, $u_G(n_4)=16r$ and $u_G(n_5)=4rq$. By Theorem \ref{1.1}, the following relationships hold for the conjugacy class size:
\begin{enumerate}
	\item [$\bullet$]
	$n_2=rq$,
	\item [$\bullet$]
	 $n_3 \in \{2q,4q,8q,16q\}$,
		\item [$\bullet$]
		 $n_4 \in \{2r,4r,8r,16r\}$,
			\item [$\bullet$]
			$n_5 \in \{2r,4r,2q,4q,2rq,4rq,rq\}$.
\end{enumerate}
Now, assume that in all cases of the set of conjugacy class sizes of $G$, denoted $cs(G)$, there exist $i \neq j$ such that $n_i = n_j$. Under this assumption, the conjugacy class sizes of $G$ will belong to one of the following: 

$\{1,rq ,2q,2r,2r\}$,
$\{1,rq ,2q,2r,2q\}$,
$\{1,rq ,2q,4r,4r\}$,
$\{1,rq ,2q,4r,2q\}$,
$\{1,rq ,2q,4r,rq\}$,
$\{1,rq ,2q,8r,2q\}$,
$\{1,rq ,2q,8r,rq\}$,
$\{1,rq ,2q,16r,2q\}$,
$\{1,rq ,2q,16r,rq\}$,
$\{1,rq ,4q,2r,2r\}$,
$\{1,rq ,4q,2r,4q\}$,
$\{1,rq ,4q,2r,rq\}$,
$\{1,rq ,4q,4r,4r\}$,
$\{1,rq ,4q,4r,4q\}$,
$\{1,rq ,4q,4r,rq\}$,
$\{1,rq ,4q,8r,4q\}$,
$\{1,rq ,4q,8r,rq\}$,
$\{1,rq ,4q,16r,4q\}$,
$\{1,rq ,4q,16r,rq\}$,
$\{1,rq ,8q,2r,2r\}$,
$\{1,rq ,8q,2r,rq\}$,
$\{1,rq ,8q,4r,4r\}$,
$\{1,rq ,8q,4r,rq\}$,
$\{1,rq ,8q,8r,rq\}$,
$\{1,rq ,8q,16r,rq\}$,
$\{1,rq ,16q,2r,2r\}$,
$\{1,rq ,16q,2r,rq\}$,
$\{1,rq ,16q,4r,4r\}$,
$\{1,rq ,16q,4r,rq\}$,
$\{1,rq ,16q,8r,rq\}$ and
$\{1,rq ,16q,16r,rq\}$.
As an example, we show a case where a contradiction \linebreak  occurs, with the other cases following similar argument. Let’s consider the first set, $\{1,rq ,2q,2r,2r\}$. In this case, we have $n_4=n_5=2r$. How ever, we know that there are  $16r$ elements with a conjugacy class size of $2r$, as well as $4rq$ elements with the same size. This leads to a equation:
\begin{align*}
 4rq=16r.	
\end{align*}
Solving this gives $4=q$, which is a contradiction. Since $q$ is an odd prime number. Therefore, this assumption leads to a contradiction, and the group cannot have these same-size conjugacy class sets. 
\end{proof}
\begin{lemma}\label{2.4}
There does not exist any finite simple group $G$ such that \linebreak $U(G)=\{1,rq,16q,4qr,16r\}$, where $q,r$ are distinct odd prime numbers.
\end{lemma}
\begin{proof}
Assume, for the sake of contradiction, that such a group exists. According to Lemma \ref{2.3}, the group $G$ has conjugate type rank 4. Furthermore, by Theorem \ref{1.3} and Theorem \ref{1.8}, we have that $G \cong PSL(2,l)$, where $l > 7$ and $l \neq 2$. Since the order of $G$ is given by $\mid G \mid = \sum_{i \in I}u_G(n_i)$, it follows that $G$ must have even order. 
From Theorem \ref{1.1}, we have $n_2=rq$. Moreover, by Proposition \ref{1.4}, the prime number set associated with $\pi(G)=\{2,r,q\}$. Therefore, $G$ is a simple $k_3$-group. By Theorem \ref{1.5}, $G$ must belong to one of the following groups: $PSL(2,9)$, $PSL(2,8)$, $PSL(2,7)$ and $PSL(2,17)$.
However, upon examining the same-size conjugacy class set of these groups in GAP, we find the following corresponding sets for each group: 
\begin{enumerate}
	\item [$\bullet$]
	$U(PSL(2,7))=\{1,rq,2rq,16r,8q\}$ ,
	\item [$\bullet$]
	$U(PSL(2,8))=\{1,r^{2}q,8r^{3},32q\}$,
	\item [$\bullet$]
	$U(PSL(2,9))=\{1,r^{2}q,16q,2r^{2}q,16r^{2}\}$,
	\item [$\bullet$]
	$U(PSL(2,17))=\{1,r^{2}q,32r^{2},2r^{3}q,2^{6}q\}$.
\end{enumerate}
On the other hand, the same-size conjugacy class of $G$ is given by \linebreak $U(G)=\{1,rq,16q,4qr,16r\}$. Clearly, none of the same-size conjugacy class sets of the groups $PSL(2,9)$, $PSL(2,8)$, $PSL(2,7)$ or $PSL(2,17)$ matches the same-size \linebreak conjugacy class set of $G$. Therefore, we reach a contradiction. No finite simple groups can satisfy the conditions set forth in this argument.
\end{proof}
In the following lemma, we will demonstrate that it is not possible for $p=q$ and $r$ to be a difference of $p$, or for $p=r$ and $q$ to be a difference of $p$. In these cases, the set $U(G)$  takes the form:
\begin{align*}
U(G) = \{1,rq,4rq,8rq,8r^2\}.	
\end{align*}
 Thus, we have the following result:
\begin{lem}\label{2.5}
	There does not exist any finite simple group $G$ such that \linebreak $U(G) =\{1,rq,4rq,8rq,8r^2\}$, where $r$ and $q$ are distinct odd primes.	
\end{lem}
\begin{proof}
	Let $G$ be a simple group with $U(G)=\{1,rq,4rq,8rq,8r^2\}$. From Theorem \ref{1.1}, we have that $n_2=rq$ and we know that $n_2 \mid \mid G \mid$, so $rq \mid \mid G \mid$. Since the order of $G$ is given by $\mid G \mid = \sum_{i \in I}u_G(n_i)$, it follows that the simple group $G$ must have even order, so $2 \mid \mid G \mid$. Moreover, by Proposition \ref{1.4}, the prime number set of order $G$ is $\pi(G)=\{2,r,q\}$. Therefore, $G$ is a simple $k_3$-group. By Theorem \ref{1.6}, $G$ must be belong to one of the following groups:
	\begin{align*}
A_5, A_6, PSL(2,7), PSL(2,8), PSL(2,17), PSL(3,3), U_3(3),\quad and \quad U_4(2).	
	\end{align*}
Using GAP, we determine the corresponding same-size conjugacy class  sets for these groups as follows:
	\begin{enumerate}
		\item [$\bullet$]
		$U(A_5)=\{1,15,20,24\}$,
		\item [$\bullet$]
		$U(A_6)=\{1,45,80,90,144\}$,
		\item [$\bullet$]
		$U(PSL(2,7))=\{1,21,42,48,56\}$,
		\item [$\bullet$]
		$U(PSL(2,8))=\{1,63,216,224\}$,
		\item [$\bullet$]
		$U(PSL(2,17))=\{1,153,288,918,1088\}$,
		\item [$\bullet$]
		$U(PSL(3,3))=\{1,104,117,624,936,1728,2106\}$,
		\item [$\bullet$]
		$U(U_3(3))=U(PSU(3,3))=\{1,56,189,378,672,1512,1728\}$,
		\item [$\bullet$]
		$U(U_4(2))=U(PSU(4,2))=\{1,45,80,240,270,480,540,720,1440,3240,5184,\linebreak 5760,6480\}$.
		\end{enumerate}
Upon examining these sets, we observe that none of the same-size conjugacy class sets matches $U(G) = \{1,rq,4rq,8rq,8r^2\}$. This leads to a contradiction, as no finite simple group satisfies the conditions imposed in this argument.
\end{proof} 
Now, we can conclude the following corollary:
\begin{cor}\label{2.6}
	Let $G$ be a finite simple group with $U(G)=\{1,rq,8pq,4qr,8pr\}$, where $p,q,r$ are prime numbers. Then $p,q,r$ must be distinct odd primes. This conclusion follows from the analysis of all possible cases regarding the relationships between 2, $p$, $q$ and $r$, where all cases except the one where $p$ , $q$ and $r$ are distinct odd primes lead to contradictions.
\end{cor}
We claim that the simple group $G$ with $U(G)=\{1,rq,4rq,8pq,8pr\}$ must be a $k_4$-group.
\begin{lemma}\label{2.7}
	Let $G$ be a simple group with $U(G)=\{1,rq,8pq,4qr,8pr\}$, where $p,q,r$ are prime numbers. Then $G$ is a $k_4$-group.
\end{lemma}
\begin{proof}
	Assume that $u_G(n_1)=1$, $u_G(n_2)=rq$, $u_G(n_3)=8pq$, $u_G(n_4)=4qr$ and $u_G(n_5)=8pr$. By Corollary \ref{2.6}, the numbers $p$, $q$ and $r$ are distinct odd primes. Since $G$ is a simple group, by Theorem \ref{1.1} and Proposition \ref{1.2}, we know that $n_2=rq$. The order of $G$ is given by $\mid G \mid = \sum_{i \in I}u_G(n_i)$, which implies that the simple group $G$ must have even order, so $2\mid \mid G \mid$. From Proposition \ref{1.4}, 2, $q$ and $r$ are elements of the prime set $\pi(G)$. Now, suppose that $p \nmid \mid G \mid$. 
	In this case, the simple group $G$ is a $k_3$-group. The order of the same-size conjugacy set of $G$ is 5, and by Theorem \ref{1.6}, the group $G$ must be one of the following groups:
	\begin{align*}
		A_6 , PSL(2,7) \quad  or \quad PSL(2,17).
	\end{align*}
\begin{enumerate}
	\item [$\bullet$]
	 Suppose $G \cong A_6$. Using GAP, we find the same-size conjugacy class set of $A_6$ to be:
	\begin{align*}
		U(A_6)=\{1,45,80,90,144\}.
	\end{align*}
\item[$\bullet$] 
For $PSL(2,7)$, we have:
\begin{align*}
	U(PSL(2,7))=\{1,21,42,48,56\}.
\end{align*}
\item[$\bullet$] 
For $PSL(2,17)$, the conjugacy class size set is:
\begin{align*}
	U(PSL(2,17))=\{1,153,288,918,1088\}.
\end{align*}
\end{enumerate}	
 None of the numbers in these sets is the product of two prime numbers. Therefore, $G$ can not be $A_6$ or $PSL(2,17)$, as neither satisfies the condition \linebreak $U(G)= \{1,rq,8pq,4qr,8pr\}$. Next, assume $G \cong PSL(2,7)$. In this case, we would have \linebreak$U(G)=U(PSL(2,7))$, which leads to a contradiction, as $U(G)$ can not match \linebreak $U(PSL(2,7))$ exactly. Therefore, we conclude that $p \in \pi(G)$, and thus $G$ is a \linebreak$k_4$-group.	 
\end{proof}
\begin{lemma}\label{2.10}
	All the simple $k_4$-groups with $\mid U(G) \mid = 5$ are isomorphic to $PSL(2,q)$, where $q$ is a prime power satisfying:
	
	$q(q^2 - 1) = gcd(2,q-1)2^{\alpha_1}3^{\alpha_2}s^{\alpha_3}t^{\alpha_4}$, with $s>3$ and $t>3$ distinc prime number.
\end{lemma}
\begin{proof}
	By using GAP, and calculating the size of $U(G)$ for all of the groups in Theorem \ref{1.10}, we have all groups except $PSL$ have the size of $U(G)$ not equal to 5 .
\end{proof}
\begin{theorem}\label{2.9}
Let $G$ be a finite simple group with $U(G)=\{1,rq,8pq,4qr,8pr\}$, where $p,q,r$ are prime numbers. Then $G \cong PSL(2,l)$, where $l$ is an odd prime greater than 7 and is uniquely determined.
\end{theorem}
\begin{proof}
By Corollary \ref{2.6} and Lemma \ref{2.7}, the simple group $G$ is a $k_4$-group. From \linebreak Lemma \ref{2.10}, $G$ is isomorphic to $PSL(2,l)$, where $l$ is an odd prime power greater than 5. By Theorem \ref{1.8}, $l$ must be a prime number greater than 7. Now, suppose $G \cong PSL(2,l_1)$ and $G \cong PSL(2,l_2)$ both share the same-size conjugacy class set. Then,the centralizers satisfy:
\begin{align*}
	\mid Cent(G) \mid = \mid Cent(PSL(2,l_1)) \mid = \mid Cent(PSL(2,l_2)) \mid	
\end{align*}
By Remark 2.2 from \cite{Zarrin2}, this implies $l_1 = l_2$. Hence, the value of $l$ must be unique. 
\end{proof}
\begin{theorem}\label{2.8}
	Let $G$ be a finite simple group with $U(G)= U(PSL(2,11))=\linebreak \{1,55,120,220,264\}$. Then $G \cong PSL(2,11)$. 
\end{theorem}
\begin{proof}
	Assume $r=11$, $q=5$, and $p=3$. Consider the set $\{1,rq,8pq,4qr,8pr\}$, which evaluates to $\{1,55,120,220,264\}$. Therefore $U(PSL(2,11))=U(G)=\linebreak \{1,rq,8pq,4qr,8pr\}$. 
	By Theorem \ref{2.9}, group $G$ is isomorphic to $PSL(2,l)$, where $l$ is an odd prime number greater than 7. So $\mid G \mid = \mid PSL(2,l) \mid$. 
	Given $U(G) = U(PSL(2,11))$ and we know $\sum u_G(n_i) = \mid G \mid$, it follows that $\mid G \mid = \mid PSL(2,11) \mid = 660$. 
	 We also know that $\mid PSL(2,l) \mid = l(l^2 - 1) / 2$. Solving $l(l^2 - 1) / 2 = 660$, yields $l=11$. Therefore $G \cong PSL(2,11)$. 
\end{proof}

\end{document}